\documentclass[12pt]{amsart}

\newcommand{\area}{{\rm Area}}

\usepackage{amsthm,amsmath,amsfonts, amssymb}
\usepackage[all]{xy}
\usepackage{psfrag}
\usepackage[dvips]{graphicx}

\usepackage[latin1]{inputenc}

\makeatletter
\def\th@alexnormal{%
\let\thm@indent\noindent 
\thm@headfont{\bfseries}
\normalfont
}
\makeatother

\makeatletter
\def\th@alexit{%
\let\thm@indent\noindent 
\thm@headfont{\bfseries}
\normalfont
\fontshape{it}
\selectfont
}
\makeatother

\theoremstyle{alexit}

\newtheorem{theorem}[equation]{Theorem}

\theoremstyle{remark}

\theoremstyle{definition}

\newtheorem{question}[equation]{Question}
\numberwithin{equation}{section}

\begin{document}
\author{Alexandre Girouard}
\address{Institut de Math\'ematiques de Neuch\^atel, rue \'Emile-Argand 11,
  2000 Neuch\^atel, Suisse.}
\email{alexandre.girouard@unine.ch}

\author{Iosif Polterovich}
\address{D\'epartement de math\'ematiques et de
statistique, Universit\'e de Montr\'eal, C. P. 6128,
Succ. Centre-ville, Montr\'eal, Qu\'ebec,  H3C 3J7,  Canada}
\email{iossif@dms.umontreal.ca}

\title{Upper bounds for Steklov eigenvalues on surfaces}


\begin{abstract}
  We give explicit isoperimetric upper bounds for all Steklov
  eigenvalues of a compact  orientable surface with boundary, in terms of the genus, the length of the
  boundary, and the number of boundary components.  Our estimates
  generalize a recent result of Fraser--Schoen, as well as the
  classical inequalites obtained by Hersch--Payne--Schiffer, whose
  approach is used in the present paper.
\end{abstract}

\maketitle

\section{Introduction}
\subsection{Steklov spectrum}
Let $\Sigma$ be a compact orientable surface with
boundary,  and let
$\Delta$  
be the Laplace--Beltrami operator associated with a Riemannian
metric on $\Sigma$. The Steklov eigenvalue problem on $\Sigma$ is given by:
\begin{gather*}
  \Delta u =0 \,\,\, {\rm in}\,\,\,  \Sigma, \quad \partial_n u = \sigma u \,\,\, {\rm on} \,\,\, \partial \Sigma,
\end{gather*}
where  $\partial_n$ denotes the outward normal derivative. The
spectrum of the Steklov problem is discrete and its eigenvalues form
a sequence 
$$0=\sigma_0<\sigma_1\leq\sigma_2\leq\sigma_3\leq\cdots\nearrow\infty,$$
where each eigenvalue is repeated according to its multiplicity~\cite{band}.
The eigenfunctions $\phi_k$, $k=0,1,2\dots$  can be chosen to form an
orthogonal basis of $L^2(\partial \Sigma)$. Note that the eigenfunction $\phi_0$ 
corresponding to $\sigma_0=~0$ is constant.

The Steklov eigenvalues coincide with the eigenvalues of  the
Dirichlet-to-Neumann map $\Lambda$. If the boundary $\partial \Sigma$ is smooth, it is a pseudo-differential elliptic operator 
$\Lambda:C^\infty(\partial\Sigma)\rightarrow C^\infty(\partial\Sigma)$ of order one \cite{Taylor},
defined by
$$\Lambda(f)=\partial_nHf,$$
where $Hf$ is the harmonic extension of $f$ to the interior of
$\Sigma$ (i.e. $\Delta(Hf)=0$ on $\Sigma$). 
The Dirichlet-to-Neumann map has important applications to inverse
problems~\cite{cldr,LaTaUh}. 

\subsection{Main results}
Isoperimetric inequalities for Steklov eigenvalues have been actively
studied for more than fifty years
\cite{wein,WheeHorg,brock, HenPhil, gp}. In particular,  a number of recent
papers are concerned with the Steklov spectrum on manifolds with
boundary \cite{esco2,fraschoen, hassannezhad, ceg2}.  The following
estimate on the first nontrivial Steklov eigenvalue on a surface with
boundary was proved by Fraser and Schoen~\cite{fraschoen}:
\begin{equation}
  \label{fsineq}
  \sigma_1\,L\leq 2\pi(\gamma+l).
\end{equation}
Here $L$ is the length of the boundary, $\gamma$ is the genus of the surface and $l$ the number of boundary components.
For simply connected planar domains, inequality~\eqref{fsineq}  is sharp and was
proved  by Weinstock in~\cite{wein}.

The goal of this note is to generalize  \eqref{fsineq} to higher eigenvalues. We prove
\begin{theorem}\label{CoroHPSgenus} Let $\Sigma$ be a compact orientable surface of genus $\gamma$, such that the
 boundary $\partial\Sigma$ has  $l$ connected components of total length $L$. Then 
  \begin{gather}
    \label{CoroIneq}
    \sigma_k\,L\leq2\pi(\gamma+l)\,k
  \end{gather}
for any integer $k\ge 1$.
\end{theorem}
In fact, Theorem \ref{CoroHPSgenus} is a special case  (set  $p=q$
below)  of the following result:
\begin{theorem}\label{thmHPSgenus}
  Under the assumptions of Theorem \ref{CoroHPSgenus},
  \begin{gather}\label{mainIneq}
    \sigma_p\sigma_q\,L^2\leq
    \begin{cases}
      \pi^2(\gamma+l)^2 (p+q)^2 &\mbox{ if }p+q\mbox{ is even},\\
      \pi^2(\gamma+l)^2  (p+q-1)^2 &\mbox{ if }p+q\mbox{ is odd},
    \end{cases}
  \end{gather}
 for any pair of integers $p,q \ge 1$.
\end{theorem}

\subsection{Discussion} It follows from Weyl's law for eigenvalues of the Dirichlet-to-Neumann operator that 
the linear dependence on $k$  in \eqref{CoroIneq} is optimal. For simply connected planar domains, the inequalities \eqref{mainIneq}
were obtained by Hersch, Payne, and Schiffer
in~\cite{HPS}. In~\cite{gp} we proved that in this case (here $\gamma=0$, $l=1$)
the estimates~\eqref{CoroIneq} are sharp for all $k\ge 1$. 
We do not expect \eqref{CoroIneq} to be sharp for other values of $\gamma$ and $l$ (cf. \cite[Theorem 2.5]{fraschoen});
see also Question \ref{quest} below.

The proof of Theorem
\ref{thmHPSgenus}  combines the methods of  \cite{fraschoen} and
\cite{HPS}.  Following \cite{fraschoen},  we use a version of
Ahlfors Theorem~\cite{MR0036318} proved by Gabard~\cite{gabard},
according to which any Riemannian surface of genus $\gamma$ with $l$ boundary components can be
represented as a proper conformal branched cover of a disk $\mathbb{D}$ with
degree at most $\gamma+l$. Properness of the covering map implies
that the boundary $\partial\Sigma$ is mapped to the circle $S^1$.  It
is essential in this proof since the test functions for the
variational characterization of the eigenvalues $\sigma_k$ are built
from the eigenfunctions of a  certain  one--dimensional problem on
$S^1$. This approach  was suggested by Hersch, Payne and Schiffer
in~\cite{HPS}.

The analogue of the estimate \eqref{fsineq} for the first nonzero
Laplace eigenvalue $\lambda_1$ on a closed surface $\Sigma$ (without
boundary) is the Yang--Yau inequality~\cite{YY}~:
\begin{equation}
\label{YY}
\lambda_1 \,\area (\Sigma) \le 8\pi d,
\end{equation}
where $d$ is the degree of a conformal branch covering of $\Sigma$
over a sphere. It was observed in~\cite{MR1046044} that one could take
$d\leq [\frac{\gamma+3}{2}]$, where $[\cdot]$ denotes the integer part.

For higher eigenvalues of the Laplacian 
on surfaces,  no explcit estimates like \eqref{CoroIneq} are
known. However,  with an implicit constant such a bound was proved by
Korevaar in~\cite{kvr} using a different approach. The analogue of
Korevaar's result for Steklov eigenvalues on surfaces was obtained
in~\cite{ceg2} (see also~\cite[Section 5.3]{gyn} and \cite[Example 1.3]{kok}): there exists a
universal constant $C$ such that
\begin{gather}\label{ineqCEG}
  \sigma_kL\leq C\, (\gamma+1)k,\quad k=1,2,3,\dots
\end{gather}
Note that the bound \eqref{ineqCEG} does not depend on the number of
boundary components of $\partial \Sigma$, which makes it a sharper estimate than
\eqref{CoroIneq} for $l$ large enough. Another
interesting development of Korevaar's method for both Laplace and
Steklov eigenvalues can be found in \cite{hassannezhad} where
$\lambda_k$ and $\sigma_k$ are bounded by a linear combination of $k$
and $\gamma$ (instead of its product). However, the constants
in~\cite{hassannezhad} are also implicit.


Let us conclude by an open question.  It was proved 
in~\cite{Buser}
that there exists a sequence of closed surfaces $\Sigma_n$ of genera
$\gamma_n \to \infty$ such that
$$\lim_{n\rightarrow\infty}\lambda_1(\Sigma_n)\,\area(\Sigma_n)=\infty.$$
Moreover,  it was subsequently shown 
in \cite{BrooksMakover}
that one can choose a sequence  of surfaces with $\gamma_n=n$ and 
$\lambda_1(\Sigma_n)\,\area(\Sigma_n)$ growing  linearly as $n\nearrow\infty$.
Therefore, the dependence on the genus $\gamma$ in the Yang--Yau inequality~\eqref{YY}
is optimal up to a multiplicative constant.
\begin{question}
\label{quest}
  Is there a sequence $\Sigma_n$ of surfaces with boundary of genera
  $\gamma_n \to \infty$ such that 
  $\sigma_1(\Sigma_n)L(\partial\Sigma_n) \to \infty$ as $n\to\infty$?
  If yes, is it possible to achieve  linear growth?
\end{question}

\subsubsection*{Acknowledgments}
The authors would like to thank Bruno Colbois for 
fruitful discussions.  Research of I.P. was supported by NSERC, FQRNT and the Canada Research Chairs program.

\section{Proof of Theorem~\ref{thmHPSgenus}}
\subsection{Reduction to the circle}
Let $\left(\phi_k\right)_{k=0}^{\infty}\subset
L^2(\partial\Sigma)$ be a complete orthonormal system of 
eigenfunctions of the Dirichlet--to--Neumann map. It is well known that
if a function $f\in C^\infty(\Sigma)$ satisfies
\begin{gather}\label{orthoCondition}
  \int_{\partial\Sigma}f\phi_j\quad\mbox{ for }j=0,1,2\dots,k-1,
\end{gather}
then 
\begin{gather}\label{VarCharac}
  \sigma_k\leq R_\Sigma(f):=
  \frac{\int_{\Sigma}|\nabla f|^2}{\int_{\partial\Sigma}f^2}.
\end{gather}

\medskip
The proof of Theorem~\ref{thmHPSgenus} is based on the
approach of~\cite{HPS}. We construct test functions using linear
combinations of harmonic oscillators on $S^1$, extend them
harmonically to the disk and then lift  to a branched cover
representation of $\Sigma$. Using sufficiently many  harmonic
oscillators, one can ensure the existence of a linear combination
satisfying the orthogonality conditions~(\ref{orthoCondition}).

\medskip
As was shown in~\cite{gabard}, there exists a proper conformal
branched cover
$$\psi:\Sigma\rightarrow\mathbb{D}$$
of degree $d\leq\gamma+l$. Because $\psi$ is proper, it takes
the boundary $\partial\Sigma$ to the circle $S^1=\partial
\mathbb{D}$. The restriction of $\psi$ to each connected
component of $\partial\Sigma$ is a covering map of $S^1$.
Let $ds$ be the Riemannian measure on $\partial\Sigma$.
We define the push-forward measure
$d\mu=\psi_*ds$ on the circle $S^1$,
and introduce the ``mass parameter''
$$m(\theta)=\int_0^{\theta}d\mu(\theta).$$
In particular, $d\mu=m'(\theta)d\theta$ is absolutely continuous with
respect to the Lebesgue measure $d\theta$, and the length of the
boundary $\partial\Sigma$ is given by
$$L=m(2\pi)=\int_{S^1}d\mu.$$
Given a smooth periodic function $h:\mathbb{R}\rightarrow\mathbb{R}$ with period $L$, define $f:S^1\rightarrow\mathbb{R}$ by
$$f(\theta)=h(m(\theta)).$$
The function $f$ admits
a unique harmonic extension $u$ to the disk $\mathbb{D}$. Because the
disk is simply connected, this function has  a unique harmonic
conjugate $v$ normalized in such a way that 
\begin{gather}
  \label{norm}
  \int_{S^1}v\,d\mu=0.
\end{gather}
Let the functions $\alpha,\beta:\Sigma\rightarrow\mathbb{R}$ be defined by
\begin{gather*}
  \alpha=u\circ\psi\quad\mbox{ and }\quad\beta=v\circ\psi.
\end{gather*}
Recall that the map $\psi$ is a $d$-fold conformal branched covering of $\mathbb{D}$.
It follows from conformal invariance of the Dirichlet energy in two dimensions (see also \cite{YY})  that
\begin{gather}\label{energyCoverI}
  \int_\Sigma|\nabla\alpha|^2=d\int_{\mathbb{D}}|\nabla u|^2,
  \quad
  \int_\Sigma|\nabla\beta|^2=d\int_{\mathbb{D}}|\nabla v|^2.
\end{gather}
Moreover, the Cauchy--Riemann equations imply that these two
quantities are equal.
 Integration by parts gives
\begin{gather}\label{RHSequal}
  \int_{\mathbb{D}}|\nabla u|^2=\int_{\mathbb{D}}|\nabla v|^2=\int_{S^1}v\,\partial_rv.
\end{gather}
Multiplying the two equations in~(\ref{energyCoverI}) and
using~(\ref{RHSequal}),  we get
\begin{gather}
\label{aaa}
  \int_\Sigma|\nabla\alpha|^2\int_\Sigma|\nabla\beta|^2
  =d^2\left(\int_{S^1}v\partial_rv\right)^2.
\end{gather}
The Cauchy--Riemann equations also imply the pointwise equality 
$$\partial_rv=-\partial_\theta u=-f'(\theta)=-h'(m(\theta))m'(\theta).$$
Applying the Cauchy--Schwarz inequality to the measure
$d\mu=m'(\theta)d\theta$  leads to~:
\begin{multline}
\label{hhh}
  \left(\int_{S^1}v\partial_rv\right)^2=
  \left(
    \int_{S^1}v(\theta)h'(m(\theta))\,\overbrace{m'(\theta))d\theta}^{d\mu(\theta)}
  \right)^2\\
  \leq \int_{S^1}v^2(\theta)\,d\mu(\theta)\int_{S^1}h'(m(\theta))^2\,d\mu(\theta).
\end{multline}
At the same time,
\begin{gather}
\label{www}
  \int_{\partial\Sigma}\alpha^2\,d_{\partial\Sigma}
  =\int_{S^1}f^2\,d\mu\,\quad\mbox{ and }\quad
  \int_{\partial\Sigma}\beta^2\,d_{\partial\Sigma}
  =\int_{S^1}v^2\,d\mu.
\end{gather}
Estimating the product of the Rayleigh quotients
$R_\alpha:=R_\Sigma(\alpha)$ and $R_\beta:=R_\Sigma(\beta)$ using the relations \eqref{aaa}, \eqref{hhh} and \eqref{www},  
we notice that   $\int_{S^1} v^2(\theta) \,d\mu(\theta)$ cancels out on the right--hand side.   This is the key trick  in the method introduced in \cite{HPS}. Namely, we obtain the following bound:
\begin{gather*}
  R(\alpha)R(\beta)\leq d^2
  \frac{\int_{S^1}h'(m(\theta))^2\,d\mu(\theta)}
  {\int_{S^1}h(m(\theta)^2\,d\mu(\theta)}
= d^2R_L(h).
\end{gather*}
Here
$$R_L(h):=\frac{\int_{0}^Lh'(m)^2\,dm}
  {\int_{0}^Lh(m)^2\,dm}$$
is the Rayleigh quotient of a uniform circular string of length
$L$. Its eigenmodes are well known.
Let  $h_k:~\mathbb{R}\rightarrow\mathbb{R}$, $k=0,1,2\dots$,  be defined by
$h_0=1$ and
\begin{gather*}
  h_k(m)=
  \begin{cases}
    \cos\left(\frac{2n\pi m}{L}\right)&\mbox{ if }k=2n-1,\\
    \sin\left(\frac{2n\pi m)}{L}\right)&\mbox{ if }k=2n.
  \end{cases}
\end{gather*}
 for $k\geq 1$.
Clearly,
\begin{gather*}
  R_L(h_k)=\left(\frac{2\pi n}{L}\right)^2\quad\mbox{for }k=2n\mbox{
    or }k=2n-1.
\end{gather*}
This leads to
\begin{gather}\label{ProductRayleigh}
  R(\alpha)R(\beta)\leq
  \left(\frac{\pi d}{L}\right)^2
  \begin{cases}
    k^2&\mbox{ if }k=2n,\\
    (k+1)^2&\mbox{ if }k=2n-1.
  \end{cases}
\end{gather}

\subsection{Construction of test-functions}
The rest of the argument is almost exactly the same as
in~\cite{HPS}. We present it  below  for the sake of completeness.
Let $N=p+q-1$. Consider a function
\begin{gather}\label{linComb}
  f=\sum_{k=1}^{N}c_kf_k,\quad (c_k\in\mathbb{R}),
\end{gather}
where the functions $f_k: \mathbb{S}^1 \to \mathbb{R}$ are defined by $f_k(\theta)=h_k(m(\theta)$. 
The functions $f_k$ are $d\mu$--orthogonal to each other, and hence
linearly independent.  The harmonic extensions $u_k$ of $f_k$ are also
linearly independent, because taking the harmonic extension is a linear and
injective operation. For the same reason, the harmonic conjugates
$v_k$  are linearly independent as well. Moreover, since by definition
$f_0=1$,  $f_k$ are $d\mu$--orthogonal to constants for all
$k=1,2,3,\dots$,  and hence $\int_{\partial \Sigma} \alpha_k=0$ for all
$k\ge 1$,  where $\alpha_k=u_k \circ \psi$.
At the same time, by the normalization \eqref{norm}, 
$\int_{\partial \Sigma} \beta_k=0$ for all $k\ge 1$,
where  $\beta_k=v_k \circ \psi$.
Let 
\begin{gather*}
  u=\sum_{k=1}^{N}c_ku_k\quad\mbox{ and }\quad v=\sum_{k=1}^{N}c_kv_k.
\end{gather*}
As before, these functions are lifted to $\alpha=u\circ\psi$ and
$\beta=v\circ\psi$.

In order to use $u$ and $v$ in the variational
characterization~(\ref{VarCharac}) for $\sigma_p$ and $\sigma_q$
respectiveley, they have to satisfy the orthogonality
conditions~(\ref{orthoCondition})~:
\begin{gather*}
  \int_{\partial\Sigma}\alpha\phi_k=0\quad\mbox{ for } k=1,\cdots,p-1\\
  \int_{\partial\Sigma}\beta\phi_k=0\quad\mbox{ for } k=1,\cdots,q-1
\end{gather*}
These $N-1$ linear constraints can be resolved for some choice of
$N$ constants $c_1,\dots,c_N$. It follows
from~(\ref{ProductRayleigh}) that 
\begin{gather*}
  \sigma_p\sigma_q\leq R(\alpha)R(\beta)\leq
  d^2R_L(h),
\end{gather*}
where $h=\sum_{k=1}^Nc_kh_k$.
We conclude by observing that
\begin{align*}
  R_L(h)\leq R_L(h_N)&=
  \left(\frac{\pi d}{L}\right)^2
  \begin{cases}
    N^2&\mbox{ if }N \mbox{ is even},\\
    (N+1)^2&\mbox{ if }N \mbox{ is odd}.
  \end{cases}\\
  &=
  \left(\frac{\pi d}{L}\right)^2
  \begin{cases}
    (p+q-1)^2&\mbox{ if }p+q \mbox{ is odd},\\
    (p+q)^2&\mbox{ if }p+q \mbox{ is even}.
  \end{cases}
\end{align*}
Recalling that $d\leq\gamma+l$ completes the proof of
Theorem~\ref{thmHPSgenus}. \qed
\bibliographystyle{plain}
\bibliography{biblioHPS4}
\end{document}